\documentclass{amsart}
\usepackage{hyperref}

\newtheorem{theorem}{Theorem}[section]
\newtheorem{lemma}[theorem]{Lemma}

\newtheorem{remark}[theorem]{Corollary}
\numberwithin{equation}{section}

\begin{document}

\title[Complex perturbations of periodic matrices]{
The discrete spectrum for complex
perturbations of periodic Jacobi matrices}

\author{Iryna Egorova}
\address{Institute for Low Temperature Physics and Engineering\\ 47\\ Lenin Ave\\
Kharkiv \\ Ukraine}
\email{\href{mailto:egorova@ilt.kharkov.ua}{egorova@ilt.kharkov.ua}}

\author{Leonid Golinskii}
\address{Institute for Low Temperature Physics and Engineering\\ 47\\ Lenin Ave\\
Kharkiv \\ Ukraine}
\email{\href{mailto:golinskii@ilt.kharkov.ua}{golinskii@ilt.kharkov.ua}}

\thanks{{\it The work was partially supported by INTAS grant no. 03-51-6637.}}


\begin{abstract} We study spectrum inclusion regions for
complex Jacobi matrices which are compact perturbations of real periodic
Jacobi matrix. The condition sufficient for the lack of discrete spectrum
for such matrices is given.
\end{abstract}

\maketitle

 \section{Introduction}

\indent

The main object under investigation is an infinite complex Jacobi
matrix $J$, which is a compact perturbation of a real symmetric
periodic matrix $J_0$ of period $N$ (background). The spectrum of
the background as an operator in $\ell^2(\mathbb{N})$ is known to
be a union of at most $N$ bands divided by gaps, each of which may
carry at most one eigenvalue. Since compact perturbations preserve
the continuous spectrum we are primarily interested in the
location of the discrete spectrum of $J$. According to the general
results of operator theory for small in norm perturbations the
spectrum of perturbed operator is settled in a small vicinity of
the spectrum of the background. A specific feature of the given
situation is that the spectrum is accumulated near the  edges of
bands and near the  discrete spectrum of the background. Whereas a
qualitative picture holds for any value of $N$, we  find certain
quantitative bounds for the size of domain where the discrete
spectrum lives only for $N=2$. In this case we also show that the
discrete spectrum of $J$ is empty provided the spectrum of the
background is pure continuous and the first moment of perturbation
is small enough.

We proceed as follows. In Section 2 we present some basic facts
from the spectral theory of periodic Jacobi matrices. The key tool
of our investigation - the discrete integral equation for the Jost
solution of the  Jacobi difference equation (three term recurrence
relation) for the perturbed matrix $J$ - is introduced in Section
3. In Section 4 we describe the algorithm of reconstruction of a
2-periodic matrix from its Borg sequence (lemma \ref{theor4.1})
and obtain a new representation for the Floquet solution (lemma
\ref{theor4.2}). Finally, in the last Section 5 we prove the main
results of the paper.

Let
$$
    J=\left(\begin{array}{ccccc}
                b_1&c_1&&&\\
                a_1&b_2&c_2&&\\
                &a_2&b_3&c_3&\\
                &&\ddots&\ddots&\ddots\\
            \end{array}\right), \qquad
    J_0=\left(\begin{array}{ccccc}
                b_1^0&a_1^0&&&\\
                a_1^0&b_2^0&a_2^0&&\\
                &a_2^0&b_3^0&a_3^0&\\
                &&\ddots&\ddots&\ddots\\
            \end{array}\right),            $$
be Jacobi matrices with the  entries
 $$ a_n,\ b_n,\ c_n\in\mathbb{C},\quad a_nc_n\neq 0; \qquad
 a_{n+2}^0=a_n^0>0,\quad b_{n+2}^0=b_n^0\in\mathbb{R},$$
and
\begin{eqnarray*}
  (J y)_n &=& a_{n-1}y_{n-1} +b_n y_n + c_n y_{n+1} = \lambda y_n,
  \quad \lambda\in\mathbb{C} \\
  (J_0 \,y)_n &=& a_{n-1}^0\,y_{n-1} +b_n^0\,y_n + a_n^0\,y_{n+1} = \lambda y_n,
 \end{eqnarray*}
be the  Jacobi difference equations, which correspond to the
matrices $J$ and $J_0$, respectively. We assume that either the
zero moment $\sum_m d_m$ or the first moment $\sum_m md_m$ of the
perturbation is finite, where
\begin{equation}\label{0.1}d_m:= |b_m^0 - b_m| + \left|a_{m-1}^0 -
\frac{a_{m-1}c_{m-1}}{a_{m-1}^0}\right|.\end{equation}

With no loss of generality we may suppose that the continuous spectrum of
the background $J_0$ is a set of the form
$$\sigma_c(J_0)=[-d,-s]\,\bigcup \,[s,d],\quad 0<s<d.$$

Put $\varepsilon=1$ when $J_0$ has an eigenvalue $\nu$, $ -s< \nu
<s$, and $\varepsilon=-1$ otherwise (this is actually the case
when $\nu$ is an eigenvalue of the Jacobi matrix with shifted
entries $\hat a_n^0=a_{n+1}^0$, $\hat b_n^0=b_{n+1}^0$ or $\nu=\pm
s$)\footnote{See Section 2 for more information about periodic
problem}. Let $w^2(\lambda)$, $\, |w(\lambda)|\leq 1$, be the
Floquet multiplier (see (\ref{1.6})) of $J_0$. The following
result holds.

\begin{theorem}\label{theor1} Let $\sum_{m=1}^\infty d_m<\infty$, where the value
$d_m$ is given in $(\ref{0.1})$.
 Then in the domain of variable $\lambda$, defined by the relation
\begin{equation}
\frac{K_2(s,d)|w(\lambda)|}{|w^2(\lambda) -
w^{2\varepsilon}(\nu)||1-w^4(\lambda)|}\,\sum_{m=1}^\infty\,d_m<t,\label{0.2}
\end{equation}
with $t\simeq 0.567..$  being a unique positive root of the
equation $x\exp(x)=1$ and $$K_2(s,d)=\frac{1024\sqrt
2\,d^5}{(d-s)^5\,\left(d^2 - s^2\right)^{1/2}},$$ the operator $J$
has no discrete spectrum at all.\end{theorem}

It can be seen from (\ref{0.2}), that if perturbation $\sum d_m$
is small enough, the discrete spectrum of $J$ is concentrated near
the edges of bands (which are exactly the points satisfying
equation $w^4=1$) and in a vicinity of the point $\nu$, as long as
it is an eigenvalue of $J_0$ (i.e., $\varepsilon=1$), or $\nu=\pm
s$. We can also observe that even if $\nu$ is not an eigenvalue of
$J_0$, its position relative to the edge of the gap exerts certain
influence on the size of the domain where the discrete spectrum of
$J$ is contained.

In the case when the first moment is small enough we have
\begin{theorem} Let the background $J_0$ has pure continuous spectrum and
$a_n^0\neq a_{n+1}^0$. Let $K_2(s,d)$ and $t$ be taken from theorem 1.1.
Then under assumption
$$\frac{K_2(s,d)}{1-|w(\nu)|^2}\sum_{m=1}^\infty md_m <t,$$
the complex Jacobi matrix $J$ has no discrete spectrum.
\end{theorem}

Our work was obviously very strongly influenced by Geronimo and
Van Assche \cite{GV}, wherein the behavior of the Jost solutions
for real perturbations of a periodic Jacobi matrix was
investigated (see also \cite[Section 7]{T}). More detailed
analysis of such solutions in the case $N=2$ is based upon the
reconstruction algorithm for a periodic Jacobi matrix from its
Borg sequence, suggested in \cite{Per}.

In conclusion let us point out that the similar problems for the
case of constant background ($a_n^0=1$, $b_n^0=0$) are discussed
in \cite{EG} (see also \cite{L}).

\section{Preliminaries on periodic spectral problem}

Let $$
    J_0=\left(\begin{array}{ccccc}
                b_1^0&a_1^0&&&\\
                a_1^0&b_2^0&a_2^0&&\\
                &a_2^0&b_3^0&a_3^0&\\
                &&\ddots&\ddots&\ddots\\
            \end{array}\right),\quad a_j^0>0,\quad b_j^0=\overline{b_j^0}$$
be an infinite Jacobi matrix with $N$ - periodic entries, $N\geq2$,
$$ a_{j+N}^0= a_j^0; \quad b_{j+N}^0=b_j^0; \quad j\in \mathbb{N}=\{1,2,\ldots\}. $$
Put $a_0^0 = a_N^0$ and consider the three term equation
\begin{equation} a_{n-1}^0 y_{n-1} +b_n^0 y_n + a_n^0 y_{n+1} = \lambda
y_n;\quad n\in\mathbb{N},\ \ \lambda\in\mathbb{C}.
\label{1.1}\end{equation} The solutions $s_n$ and $c_n$ of `$\sin$' and
`$\cos$'\ type, defined by the initial data $s_0=c_1=0$, $c_1=s_0=1$, are
polynomials with real coefficients of degree $n-1$ and $n-2$ respectively,
and also
$$s_n(\lambda) = \frac{\lambda^{n-1}}{a_1^0\ldots a_{n-1}^0} +\ldots\,.$$

The Wronskian $\langle s,c\rangle(n)=a_{n-1}^0(s_{n-1}c_n - s_n
c_{n-1})$ of these two solutions does not depend on $n$ and by
periodicity of $a_n^0$ $$\langle
s,c\rangle(N+1)=a_{N}^0(s_{N}c_{N+1} - s_{N+1} c_{N})=\langle
s,c\rangle(1) = -a_N^0$$ or
\begin{equation}\label{1.2}s_{N+1}(\lambda)c_{N}(\lambda) - s_{N}(\lambda)
c_{N+1}(\lambda)=1,\quad \lambda\in\mathbb{C}.
 \end{equation}

The following two polynomials with real coefficients
\begin{equation}\label{1.3} u(\lambda)
=\frac{1}{2}\left(s_{N+1}(\lambda) + c_N(\lambda)\right)
=\frac{\lambda^N}{2a_1^0\ldots a_N^0}+\ldots
\end{equation}

\begin{equation}
\label{1.4} v(\lambda) =\frac{1}{2}\left(s_{N+1}(\lambda) -
c_N(\lambda)\right)
\end{equation}
are of particular interest. The first one, $u$, known as the Hill
discriminant, characterizes the continuous spectrum of the
operator $J_0$ in the space $\ell^2(\mathbb{N})$. Indeed, as is
known by \cite[Chapter 4] {Td}, the spectrum $\sigma(J_0)
=\sigma_c(J_0)\bigcup\sigma_d(J_0)$ is a union of at most $N$
disjoint closed intervals (bands) which constitute the continuous
spectrum $\sigma_c(J_0)$ $$\sigma_c(J_0)=[\mu_0,\mu_N]\setminus
\bigcup_{k=1}^{M-1}\Gamma_k,\quad \Gamma_k=(\mu_k^-,\mu_k^+),\
M\leq N,$$ with at most one eigenvalue inside each gap $\Gamma_k$
(they form the discrete spectrum $\sigma_d(J_0)$). The continuous
spectrum is related to the Hill discriminant by $$\sigma_c(J_0) =
\{\lambda:\ -1\leq u(\lambda)\leq 1\} $$ (see \cite[Section
7.1]{Td}). Throughout the  paper we assume that $M=N$, i.e., an
equation $u^2=1$ has no multiple zeros. Hence $\mu_k^-<\mu_k^+$
and $$\{\lambda:\ u(\lambda)=\pm
1\}=\{\mu_0,\mu_1^-,\mu_1^+,...,\mu_{N-1}^+,\mu_N\}.$$ By
(\ref{1.3}) $u(\mu_N) = 1$ and $u\to +\infty$ as $\lambda\to
+\infty$ . The behavior of $u$ for $\lambda\to -\infty$ is
determined by the parity of  $N$.

Consider the function $$\left( u^2(\lambda)
-1\right)^{1/2}=\frac1{2a_1^0\ldots a_n^0}\,\left((\lambda -
\mu_0)(\lambda - \mu_N)\prod_{k=1}^{N-1}(\lambda -
\mu_k^-)(\lambda - \mu_k^+)\right)^{1/2},$$ where the branch of
the square root is taken positive for $\lambda>\mu_N$. This
function is real-valued inside the gaps $\Gamma_k$ in such a way
that $\mbox{sign}\left (u^2 - 1\right)^{1/2} =(-1)^{N-k+1},$
$\lambda\in\Gamma_k$. On the set $\sigma_c(J_0)$ it takes pure
imaginary values.

The  Weyl function $m$ for periodic Jacobi matrices $J_0$ is known
to be of the form (see \cite[formula\,(21.3.20)]{Ti})
$$m(\lambda) =\langle(J_0 - \lambda
I)^{-1}e_1,e_1\rangle=\frac{v(\lambda) - \left( u^2(\lambda) -
1\right)^{1/2}}{s_N(\lambda)}.$$ Denote by $\Lambda_r$ the set of
all poles of $m$, which is clearly a subset of the set
$\Lambda=\{\lambda_k\}_{k=1}^{N-1}$ of all roots of the
polynomials $s_N$ (they all are real and simple and $\lambda_k\in
[\mu_k^-, \mu_k^+]).$ We assume that all $\lambda_k$ stay strictly
inside the gaps, that is
\begin{equation}\label{1.20}\Lambda\subset\bigcup_{k=1}^{N-1}
\Gamma_k.\end{equation}

 The  Weyl solution $\psi_n$ of (\ref{1.1}) is defined in a standard way
 $$\psi_n(\lambda) = c_n(\lambda) + m(\lambda)s_n(\lambda);\quad
 \psi_0=1,\ \psi_1=m,$$
 $\{\psi_n\}\in \ell^2(\mathbb{N})$ for
 $\lambda\notin\sigma_c(J_0)\bigcup\Lambda_r.$

Along with the function $m$ consider the function
 \footnote{The function $\hat m$ appears to be the Weyl function of the
operator $\hat J_0:\ \ell^2(\mathbb{Z}_-)\to \ell^2(\mathbb{Z}_-)$
with periodically extended entries $a_i$, $b_i$. }
 $$\hat m(\lambda):= \frac{v(\lambda) + \left(u^2(\lambda) -
 1\right)^{1/2}}{s_N(\lambda)}.$$
It is easy to see from (\ref{1.2})-(\ref{1.4}) that
 \begin{equation}
 \label{1.5}
 v^2(\lambda_k)=u^2(\lambda_k) - 1,\end{equation}
and so one of the functions $m, \hat m$ has a pole at each point
$\lambda_k$ and the other - a removable singularity .

Let us introduce the second solution of (\ref{1.1})
 $$ \hat \psi_n(\lambda)=c_n(\lambda) + \hat
 m(\lambda)s_n(\lambda),\quad \lambda\notin\Lambda\setminus\Lambda_r.
 $$
 As the Wronskian
 $$\langle\psi,\hat\psi\rangle=a_N^0(\psi_0\hat\psi_1 -
 \psi_1\hat\psi_0)=2a_N^0\frac{\left(u^2(\lambda) - 1\right)^{1/2}}{s_N}\neq 0,$$
 the solutions $\psi$  and $\hat\psi$ are linearly independent
off the endpoints of $\Gamma_k$. Note also that
 $$\hat\psi_n(\lambda+i0)=\overline{\psi_n(\lambda+i0)}=\psi_n(\lambda-i0)$$
for $\lambda\in\sigma_c(J_0)$, and both solutions take real values
for
 $\lambda\in\mathbb{R}\setminus\sigma(J_0)$. It follows from the properties of
 the Weyl functions  $m,\ \hat m$ that for $n\neq kN,$
$k\in\mathbb{Z}_+=\{0,1,\ldots\}$ the only one of $\psi_n,\,\hat\psi_n$
may have a pole at $\lambda_k$. For $n=kN $ they both are continuous at
$\lambda_k$.

The typical feature of the  periodic Jacobi difference equation
is the existence of a special solution of (\ref{1.1}) known as the
Floquet solution $\phi_n(\lambda)$, which satisfies
\begin{equation}\label{1.10}\phi_{n+N}(\lambda) = \rho(\lambda)\phi_n(\lambda),\quad
n\in\mathbb{Z}_+.\end{equation} Here
$$\rho(\lambda)=u(\lambda)-\left(u^2(\lambda) - 1\right)^{1/2}$$
is the  Floquet multiplier (see \cite{T}), with
$|\rho(\lambda)|=1$ for $\lambda\in\sigma_c(J_0)$ and
$|\rho(\lambda)|<1$ for $\lambda\in\sigma_c(J_0)$. It is clear
that $\{\phi_n\}\in\ell^2(\mathbb{Z}_+)$ for
$\lambda\in\mathbb{C}\setminus\sigma(J_0)$, and so it agrees up to
a constant multiplier with the Weyl solution. We can put
$\phi_n(\lambda)\equiv\psi_n(\lambda)$. In addition by
(\cite[formula\,(7.19)]{T}) for the second solution
\begin{equation}\label{1.11}\hat\psi_{n+N}(\lambda)=
\rho^{-1}(\lambda)\hat\psi_n(\lambda)
\end{equation}
holds.

Following \cite{GV} we introduce the function $w$:
\begin{equation}\label{1.6}
w^N(\lambda)=\rho(\lambda).\end{equation} It defines a conformal mapping
of the upper $\mathbb{C}_+$ and lower $\mathbb{C}_-$ half-planes onto the
lower $\mathbb{W}_-$ and upper $\mathbb{W}_+$ half-disks with cuts
   \begin{equation}
   \label{1.8}
   \mathbb{W}_{\pm}=\{w: \ |w|<1;\ \pm \mbox{Im}\,
   w>0\}\setminus\bigcup_{k=1}^{N-1}\beta_k^{\pm},
   \end{equation} where
   $$\beta_k^{\pm}=\left[\exp\left(\pm i\pi(N-k)N^{-1}\right),\,\exp\left(\pm
   i\pi(N-k)N^{-1}\right)\left|w(\gamma_k)\right|\right],$$
and $\{\gamma_k\}_{k=1}^{N-1}$ are the points of local extremum of
the Hill discriminant $u$, $\gamma_k\in(\mu_k^-,\ \mu_k^+)$ (see
\cite[p.\,110]{Per}). Moreover, the upper side of the interval
$[\mu_k^+,\ \mu_{k+1}^-]$ goes to the circular arc with endpoints
$\exp\,\frac{i\pi(k-N)}{N}$ and $\exp\,\frac{i\pi(k+1-N)}{N},$ the
lower side goes to the symmetric arc on the upper semicircle. The
gap $\overline{\Gamma_k}$ viewed as a part of the boundary of
$\mathbb{C}_+$ is mapped on both sides of the cut $\beta_k^-$, in
such a way that $\gamma_k$ corresponds to the vertex of the cut
and $\mu_k^-$, $\mu_k^+$ to
$\exp\left(i\pi\left(k-N\right)N^{-1}\mp i0\right)$, respectively.
As a part of the boundary of $\mathbb{C}_+$ the gap
$\overline{\Gamma_k}$ goes to the symmetric cut $\beta_k^+$. The
infinite gaps $(-\infty,\ \mu_0)$ and $(\mu_N,\ +\infty)$ are
mapped onto $(-1,0)$ and $(0,1)$, respectively; $w(\lambda)\to 0$
for $\lambda\to\infty$. Hence the points $\mu_0,\mu_1^-,...,\mu_N$
correspond to $w^{2N}=1$ and each point from $\Lambda$ has two
symmetric images on symmetric sides of the cuts.

Consider the functions $\psi_n(\lambda)$ and $\hat\psi_n(\lambda)$
as functions of variable $w$  for
$w\in\mathbb{W}_+\bigcup\mathbb{W}_-\bigcup(-1,0)\bigcup(0,1)$.
Then they  have continuous boundary values  with
     \begin{equation}\label{1.9}\psi_n(
     \overline{w})=\overline{\psi_n(w)},\qquad  \hat\psi_n(
     \bar{w})=\overline{\hat\psi_n(w)},\quad |w|=1,
     \end{equation}
     $$\psi_n(\bar w)=\psi_n(w)\in\mathbb{R}\qquad \hat\psi(\bar
     w)=\overline{\hat\psi_n(w)}\in\mathbb{R},\quad  w\in\beta_k^{\pm}.$$
Note that by (\ref{1.10}),(\ref{1.11}) and (\ref{1.6}) the equalities
\begin{equation}\label{1.12}
\psi_n(w)=w^n\chi_n(w),\quad\hat\psi_n(w)=w^{-n}\hat\chi(w)\end{equation}
hold, where $\chi_n(w)=\chi_{N+n}(w)$, $\hat\chi_n(w)=\hat\chi_{N+n}(w)$
are $N$ - periodic functions.

\section{Integral equation }

Take an infinite Jacobi matrix
$$
    J = \left(\begin{array}{ccccc}
                b_1&c_1&&&\\
                a_1&b_2&c_2&&\\
                &a_2&b_3&c_3&\\
                &&\ddots&\ddots&\ddots\\
            \end{array}\right),\quad a_j,\ b_j,\ c_j\in\mathbb{C},\quad
            a_jc_j\neq 0,$$
and the related three term recurrence relation
\begin{equation}
\label{2.2} a_{n-1}y_{n-1} +b_ny_n + c_ny_{n+1} = \lambda
y_n;\quad \lambda\in\mathbb{C},\ n\in\mathbb{N}.
\end{equation}
We assume that the matrix entries $\{a_n\}$, $\{b_n\}$
 and $\{c_n\}$ are close to $N$ - periodic sequences $\{a_n^0\}$ è
 $\{b_n^0\}$ in the following sense
 \begin{equation}
 \label{2.3}
 \sum_{n=1}^\infty\,\{|a_n - a_n^0| + |b_n - b_n^0| + |c_n -
 a_n^0|\}<\infty.
 \end{equation}
It is convenient  to modify (\ref{2.2}) by changing the variables
$$v_m=y_m \prod_{j=m}^\infty\frac{a_j}{a_j^0},\quad m=1,2,\ldots
$$ (the infinite product converges due to (\ref{2.3})). For $v_n$
we have
 \begin{equation}\label{2.4}
 a_{n-1}^0v_{n-1} + b_nv_n +\frac{c_na_n}{a_n^0}v_{n+1}=\lambda
 v_n.
 \end{equation}
Under the Green function (for the unperturbed problem) we
mean
\begin{equation}
 \label{2.1}
 G(\lambda;n,m):=\left\{\begin{array}{ccc}
 \frac{\psi_n(\lambda)\hat\psi_m(\lambda) -
 \psi_m(\lambda)\hat\psi_n(\lambda)}{\langle\psi,\hat\psi\rangle}
 &; & m>n\\
 0& ; & m\leq n\end{array} \right.,\quad m,n\geq 0,
 \end{equation}
where the Weyl solutions $\psi_n$, $\hat\psi_n$ are defined in Section 2.
It is easy to check that
 $$a_{n-1}^0 G(\lambda;n-1,m) + (b_n^0 - \lambda) G(\lambda;n,m) +
 a_n^0 G(\lambda;n+1,m) = \delta_{n,m},$$
where $\delta_{n,m}$ is the Kronecker symbol.

 The key role throughout the rest of the paper is played by the discrete
 integral equation (cf. \cite[formulae\,(7.87)-(7.88)]{T}
 \begin{equation}
 \label{2.5}
 v_n(\lambda) = \psi_n(\lambda) + \sum_{m=n+1}^\infty\
 A(\lambda;n,m)v_n(\lambda),\quad \lambda\notin\Lambda_r.
 \end{equation}
 with the kernel
 \begin{equation}
 \label{2.51}
 A(\lambda;n,m)=(b_m^0 - b_m)G(\lambda;n,m) + \left(a_{m-1}^0 -
 \frac{a_{m-1}c_{m-1}}{a_{m-1}^0}\right)G(\lambda;n,m-1).
 \end{equation}
Let us show that each solution of (\ref{2.5}) satisfies (\ref{2.4}) and
has the following asymptotic behavior for $n\to\infty$:
$$ v_n - \psi_n\to 0. $$
Indeed, by the definition
 \begin{eqnarray*} A(\lambda;n-1,m) & = &\frac{b_n^0 -
 b_n}{a_{n-1}^0},\\
 A(\lambda;n-1,n+1) & =  & \frac{(b_{n+1}^0 - b_{n+1})(\lambda -
 b_n^0)}{a_{n-1}^0a_n^0} + \left(a_n^0 -
 \frac{a_nc_n}{a_n^0}\right).\end{eqnarray*}
 We apply the `unperturbed' operator to both parts of (\ref{2.5})
\begin{eqnarray*} (J_0 - \lambda)v_n & = & a_{n-1}^0v_{n-1} + (b_n^0 - \lambda)
v_n +
 a_n^0v_{n+1}\\  \, & = & (J_0 - \lambda)\psi_n(\lambda)
 + a_{n-1}^0\{A(\lambda;n-1,n)v_n +
 A(\lambda;n-1,n+1)v_{n+1}\}\\
 \, & + &(b_n^0-\lambda)A(\lambda;n,n+1)
 v_{n+1}+ \sum_{m=n+2}^\infty\ (J_0
 -\lambda)_n A(\lambda;n,m)v_m\end{eqnarray*} or $$a_{n-1}^0v_{n-1} + (b_n^0
 -\lambda)v_n + a_n^0 v_{n+1} = (b_n^0 - b_n)v_n + \left(a_n^0 -
 \frac{a_nc_n}{a_n^0}\right)v_{n+1},$$
which gives exactly (\ref{2.4}).

 Multiply (\ref{2.5}) through by $w^{-n}$, where $w$
 is defined in  (\ref{1.6}):
 \begin{equation}\label{2.7}
 \tilde v_n:=v_nw^{-n}=\psi_n(\lambda)w^{-n} +
 \sum_{m=n+1}^\infty\ \tilde A(\lambda;n,m)\tilde v_m,\end{equation} $$ \tilde
 A(\lambda;n,m):=A(\lambda;n,m)w^{m-n},$$
 or
 \begin{equation}\label{2.6}
 V_n(\lambda) = \sum_{m=n+1}^\infty\ \tilde
 A(\lambda;n,m)\psi_m(\lambda)w^{-m} + \sum_{m=n+1}^\infty\ \tilde
 A(\lambda;n,m)V_m(\lambda),\end{equation}
 with
 $$V_n(\lambda) = (v_n(\lambda) - \psi_n(\lambda))w^{-n}.$$

 Equation (\ref{2.6}) is analyzed thoroughly in \cite[Theorem 3]{GV},
 for real symmetric Jacobi matrices $(c_n=a_n>0,$ $b_n=\overline{b_n})$. The
 existence of the solution $V$ such that
  $$|V_n(\lambda) s_N(\lambda)|=O(1),\quad n\to\infty$$
 uniformly on compact subsets of the domain
 $\mathbb{C}\setminus\{\mu_k^{\pm}\}$ has been established
(the reason for the factor $s_N$ to enter the asymptotic relation is to
`suppress' poles of $\psi_n$).

We will be interested in the bounds of the form
$$|v_n(\lambda) - \psi_n(\lambda)|\leq\frac{ K(\lambda)w^n\kappa_0(n)}
{|1-w^{2N}|},\qquad \lambda\notin \Lambda_r,$$ wherein a
positive function $K(\lambda)$ is given explicitly in terms
of the spectral data of the background and uniformly
bounded off certain small neighborhoods of the set
$\Lambda_r$,
\begin{equation}\label{2.8}
\kappa_0(n): = \sum_{m=n+1}^\infty d_m,\quad d_m:= |b_m^0 - b_m| +
\left|a_{m-1}^0 - \frac{a_{m-1}c_{m-1}}{a_{m-1}^0}\right|.
\end{equation}
We  solve this problem completely for the case $N=2$, that is
performed later on. In the rest of this section we  show how the
method of successive approximations works in connection with
equation (\ref{2.6}). Appropriate bounds (see lemma \ref{lem5.2})
for the kernel $\tilde A$ lead to the uniform convergence of the
series
\begin{equation}\label{2.9} V_{n,1}(\lambda):= \sum_{m=n+1}^\infty \
\tilde A(\lambda;n,m)\psi_m(\lambda)w^{-m}\end{equation} on compact
subsets of the domain $\mathbb{C}\setminus(\bigcup\{\mu_k^\pm\} \bigcup
\Lambda_r)$, and moreover
$$|V_{n,1}(\lambda)| \leq \frac{K(\lambda)}{|1-w^{2N}|} \kappa_0(n).$$
Put
$$V_{n,j+1}:= \sum_{m=n+1}^\infty \tilde A(\lambda;n,m)V_{m,j}(\lambda)$$
and prove by induction that
$$|V_{n,j}(\lambda)|\leq\frac{1}{(j-1)!}\left(\frac{K(\lambda)\kappa_0(n)}
{|1-w^{2N}|}\right)^j.$$ Then the required solution of (\ref{2.6}) is
given by
\begin{equation}\label{2.10}
V_n(\lambda) = \sum_{j=1}^\infty V_{n,j}(\lambda),\qquad
|V_n(\lambda)|\leq \frac{K(\lambda)\kappa_0(n)}{|1-w^{2N}|}
\exp\left(\frac{K(\lambda)\kappa_0(n)}{|1-w^{2N}|}\right).\end{equation}
In particular for $n=0$ thanks to $\psi_0=1$ we have
$$ |v_0(\lambda) - 1|\leq \frac{K_1(w)\sum_{m=1}^\infty \ d_m}{|1-w^{2N}|}
\exp\left(\frac{K_1(w)\sum_{m=1}^\infty \ d_m}{|1-w^{2N}|}\right),$$ which
makes it possible to describe the region of the variable $w$ wherein
$v_0\neq 0$. The similar argument (with slightly different bounds) holds
under the assumption
\begin{equation}\label{2.20}
\sum_{m=1}^\infty md_m<\infty. \end{equation} It remains only to point out
that the set of zeros of the function $v_0$ agrees with the set of
eigenvalues (discrete spectrum) of the operator $J$ off $\sigma_c(J_0)$.

\section{Reconstruction of matrix of period 2 from spectral data}

Let $J_0$ be a real symmetric Jacobi matrix of period 2 with the matrix
entries $a_{j+2}^0=a_j^0$, $ b_{j+2}^0=b_{j}^0$. The first few polynomials
$s_n$, $c_n$ can be easily computed from recurrence relation (\ref{1.1}):
\begin{equation}\label{4.0}\left\{\begin{array}{llll} s_0=0,& s_1=1,&
s_2=\frac{\lambda - b_1^0}{a_1^0},&s_3=\frac{(\lambda - b_1^0)(\lambda -
b_2^0) - (a_1^0)^2}{a_1^0a_2^0}\,,\\ c_0=1,& c_1=0,&
c_2=-\frac{a_2^0}{a_1^0},&
c_3=\frac{b_2^0-\lambda}{a_1^0}\,.\end{array}\right.\end{equation} The
Hill discriminant
\begin{equation}\label{4.1}
u(\lambda)=\frac{1}{2}\left(s_3(\lambda) +c_2(\lambda)\right)
= \frac{\lambda^2 - (b_1^0 + b_2^0)\lambda + b_1^0b_2^0 -
(a_2^0)^2 - (a_1^0)^2}{2a_1^0a_2^0}\end{equation} is a quadratic
parabola with the positive leading coefficient. Hence the bands of
continuous spectrum $\{\lambda:\ -1\leq u\leq 1\}$ are symmetric
with respect to the vertex of this parabola and have the same
length. By shifting the spectral parameter we can put the vertex
to the origin. Denote the endpoints of the spectrum by
$\{-d,-s,s,d\},$ $0<s<d$. Then
\begin{equation}\label{4.2} u(\lambda) - 1= C(\lambda^2 -
d^2),\quad u(\lambda) +1 = C(\lambda^2 - s^2),\quad u^2 (\lambda) - 1 =
C^2(\lambda^2 - d^2)(\lambda^2 - s^2),\end{equation} whence it follows
that $C=2(d^2 - s^2)^{-1} $ and
\begin{equation}\label{4.3}
u(\lambda) = \frac{2}{d^2 - s^2} \left(\lambda^2 - \frac{d^2 +
s^2}{2}\right).\end{equation}

Let $\nu$ be a root of the polynomial $s_2,$ $\nu = b_1^0$. By the
assumption (\ref{1.20}) we have $|\nu|<s$. Put $\varepsilon = 1$
in the case when the Weyl function $m(\lambda)$ has a pole at
$\nu$ and $\varepsilon = -1$ otherwise. Our goal here is to
restore the matrix $J_0$ from its Borg sequence $\{\pm s, \pm d,
\nu, \varepsilon\}$ and to derive a  formula for the Weyl function
as a function of variable $w$ (\ref{1.6}), which is related to $u$
by the equality
\begin{equation}\label{4.4} u(\lambda) =\frac{w^2(\lambda) +
w^{-2}(\lambda)}{2}.\end{equation}

By comparing (\ref{4.1}) and (\ref{4.3}) we come to the relations for the
matrix entries:
\begin{equation}\label{4.5} b_2^0= - b_1^0 = -\nu \end{equation} and
\begin{equation}\label{4.6} a_1^0a_2^0 =\frac{d^2 - s^2}{4},\quad
\nu^2 +(a_1^0)^2 + (a_2^0)^2 = \frac{d^2 + s^2}{2},\end{equation} so that
$(a_1^0)^2$, $(a_2^0)^2$ are roots of the quadratic equation
$$x^2 - \left(\frac{d^2 + s^2}{2} - \nu^2\right)x +
\left(\frac{d^2 - s^2}{4}\right)^2 = 0.$$
Hence
\begin{eqnarray}
  \nonumber 4(a_{\min}^0)^2 &=& d^2 +s^2 - 2\nu^2 -
  2\left(d^2 - \nu^2)(s^2 - \nu^2)\right)^{1/2}_+, \\
  \label{4.7} 4(a_{\max}^0)^2 &=&d^2 + s^2 -
  2\nu^2+2\left((d^2 - \nu^2)(s^2 - \nu^2)\right)^{1/2}_+,
\end{eqnarray}
 where $a_{\min}^0=\min\{a_1^0, a_2^0\}$, $a_{\max}^0=\max\{a_1^0, a_2^0\}$
and $(\ )^{1/2}_+$ designates the arithmetic value of the square
root.

Let us show that $\varepsilon = 1\ (\varepsilon = -1)$ is
equivalent to $a_1^0 <a_2^0$ $(a_2^0 <a_1^0)$. To this end write
\begin{eqnarray*}v(\lambda) & = & \frac{1}{2}(s_3(\lambda) - c_2(\lambda))
=\frac{\lambda^2 - \nu^2 -(a_1^0)^2 + (a_2^0)^2}{2a_1^0a_2^0}\\ \,
& = & \frac{2}{d^2 - s^2}(\lambda^2 - \nu^2 -(a_1^0)^2 +
(a_2^0)^2),\end{eqnarray*} and so
\begin{equation}\label{4.8}
v(\nu)=\frac{2}{d^2 - s^2}\left( (a_2^0)^2 -(a_1^0)^2 \right).
\end{equation}
With regard to (\ref{4.7}) and (\ref{4.2}) (see also (\ref{1.5})) the
latter gives
 $$v^2(\nu) =\frac{4}{(d^2 - s^2)^2}(d^2 - \nu^2)(s^2 -
\nu^2)=u^2(\nu) - 1,$$ and therefore $|v(\nu)| =\left(u^2(\nu) -
1\right)^{1/2}_+$. In two-band case the Weyl function admits
representation $$m(\lambda) = \frac{v(\lambda) +
\left(u^2(\lambda) - 1\right)^{1/2}_+}{s_2(\lambda)},\quad
|\lambda|<s, $$ and has a pole if and only if $|v(\nu)| = v(\nu)$,
that is equivalent to $a_1^0< a_2^0$ by (\ref{4.8}). Thus we come
to the following
\begin{lemma}\label{theor4.1} The entries of a two-periodic Jacobi matrix
can be restored from the Borg sequence
$\{-d,-s,s,d,\nu,\varepsilon\}$ by the formulae: $$ b_2^0= - b_1^0
= -\nu,$$ $$4(a_{1}^0)^2 = d^2 + s^2 -
2\nu^2-2\varepsilon\left((d^2 - \nu^2)(s^2 -
\nu^2)\right)^{1/2}_+,$$ $$4(a_{2}^0)^2 = d^2 + s^2 - 2\nu^2
+2\varepsilon\left((d^2 - \nu^2)(s^2 -
\nu^2)\right)^{1/2}_+.$$\end{lemma}
\begin{remark}\label{cor}
The following inequalities hold
\begin{equation}\label{4.25}
\frac{d-s}{2}\leq a_j^0\leq \frac{d+s}{2}.
\end{equation}
\end{remark}

Our next problem is to find effective bounds for the Weyl solution, more
precisely, for the function $\chi_n(\lambda)=\psi_n(\lambda)w^{-n}$ in
terms of the variable $w$. As $\chi_{n+2}=\chi_n$, $\chi_0=\psi_0=1$, we
only need estimate the function
$\chi_1(\lambda)=w^{-1}(\lambda)\psi_1(\lambda)$. By the definition of $w$
and in view of (\ref{4.4}) we have
\begin{equation}
\label{4.18} w^2(\nu) = u(\nu) - \left(u^2(\nu)
-1\right)^{1/2}.\end{equation} Remind that in the two-band case
$\left(u^2(\nu) - 1\right)^{1/2}<0$.
\begin{lemma}\label{theor4.2}
The following representation holds
\begin{equation}
\chi_1(\lambda) =\frac{4a_1^0\,(\lambda + \nu)\,w(\lambda)}{d^2 -
s^2}\,\frac{1}{w^2(\lambda) -
w^{2\varepsilon}(\nu)}.\label{4.19}\end{equation}
\end{lemma}
\begin{proof} The wronskian of two solutions $s_n(\lambda)$ and
$\psi_n(\lambda)$ does not depend on $n$ so $$\langle\psi,\
s\rangle = a_2^0(\psi_2s_3-s_2\psi_3) = a_0^0(\psi_0s_1 -
\psi_1s_0)=a_0^0:=a_2^0,$$ that is
\begin{equation}\label{4.15}
s_3(\lambda)\psi_2(\lambda) - s_2(\lambda)\psi_3(\lambda)=1.\end{equation}
Since (see Section 1) $\psi_2=w^2\psi_0=w^2$, $\psi_3=w^2\psi_1$, it
follows from (\ref{4.15}), (\ref{1.1}), (\ref{4.0}) and theorem
\ref{theor4.1} that
\begin{eqnarray*}
\psi_1 &=& \frac{s_3\,w^2 - 1}{w^2\, s_2}=\frac{((\lambda -
b_2^0)\,s_2 - a_1^0)\,w^2\,(a_2^0)^{-1} - 1}{w^2\,s_2} \\
   &=& \frac{\lambda + \nu}{a_2^0} - \frac{a_1^0\,w^2 +
a_2^0}{a_2^0\,w^2\,s_2}=\frac{\lambda + \nu}{a_2^0} -
\frac{a_1^0\,\psi_2 +a_0^0\,\psi_0}{a_2^0\,w^2\,s_2} \\
   &=& \frac{\lambda + \nu}{a_2^0}-\frac{(\lambda - b_1^0)\,
\psi_1}{a_2^0\,w^2\,s_2}=\frac{\lambda + \nu}{a_2^0}
-\frac{a_1^0}{a_2^0\,w^2}\,\psi_1,
\end{eqnarray*}
i.e.,
$$\psi_1\,\left(1 + \frac{a_1^0}{w^2\,a_2^0}\right)=\frac{\lambda +
\nu}{a_2^0},$$ and hence
\begin{equation}\label{4.20}
\psi_1(\lambda) =\frac{(\lambda +\nu)\,w^2}{a_2^0\,w^2 +
a_1^0}=\frac{4a_1^0\,(\lambda +\nu)\,w^2}{4a_1^0\,a_2^0\,w^2 +
4(a_1^0)^2}.\end{equation} But by (\ref{4.2}), (\ref{4.3}),
(\ref{4.6}) and theorem \ref{theor4.1} $$\frac{4(a_1^0)^2}{d^2 -
s^2}=\frac{2}{d^2 - s^2}\left(\frac{d^2 + s^2}{2} - \nu^2\right) +
\varepsilon\,\frac{2}{d^2 - s^2}\,\left((d^2 - \nu^2)(s^2 -
\nu^2)\right)^{1/2}_+,$$ that is $$\frac{4(a_1^0)^2}{d^2 -
s^2}=-u(\nu) - \varepsilon\left(u^2(\nu) - 1\right)^{1/2}.$$
Relations (\ref{4.18}), (\ref{4.4}) imply
\begin{equation}\label{4.21}\frac{4(a_1^0)^2}{d^2 - s^2}=
(w(\nu))^{2\varepsilon}.\end{equation} Next we plague (\ref{4.21}) into
(\ref{4.20}) and take into account (\ref{4.6}):
\begin{equation}\label{4.22}
\psi_1(\lambda) =\frac{4a_1^0\,(\lambda + \nu)\,w^2(\lambda)}{(d^2
- s^2)(w^2(\lambda) - w^{2\varepsilon}(\nu))}\,,\end{equation}
which proves (\ref{4.19}) since $\chi_1(\lambda) =
w^{-1}(\lambda)\psi_1(\lambda)$.\end{proof}

\section{Bounds for solutions of integral equation}

The bounds in the method of successive approximations displayed in Section
3 for the case $N=2$ is based upon explicit expressions for $w$ and the
Weyl functions $m,\hat m$. We begin with the following simple result.
\begin{lemma}\label{lem5.1} Let $0\leq |\nu| <s<d$ and $\varepsilon =
\pm1$. Put
\begin{equation}
\tau(\nu,\varepsilon):=\frac{\varepsilon\,\left((d^2 -
\nu^2)(s^2-\nu^2)\right)^{1/2}-
\nu^2}{d^2}\,,\quad\tau=\tau(0,\varepsilon)=\frac{\varepsilon
s}{d}\,, \quad \xi=\frac{d^2 +s^2}{2d^2}\,.
\label{5.1}\end{equation} For $|z|<1$ consider the functions
$$f(z,\nu,\varepsilon):=1 + \tau(\nu,\varepsilon)\,z +
\left((1-z)(1-\tau z)\right)^{1/2},$$ $$ g(z):=1-\xi z +
\left(1-z)(1-\tau z)\right)^{1/2},$$ where the branch of the
square root is taken positive for $-1<z<1$. Then
\begin{equation}\label{5.2}
|f(z,\nu,\varepsilon)|\geq\frac{(1+\tau(\nu,\varepsilon))^2}{4}\geq
\frac{(1-|\tau|)^2}{4}=\frac{(d-s)^2}{4d^2}\,,\end{equation}
\begin{equation}
\label{5.3} |g(z)|\geq\frac{|\xi^2 -
\tau^2|}{4}=\frac{(1-\tau^2)^2}{4}=\frac{(d^2 - s^2)^2}{4d^4}\,.
\end{equation}\end{lemma}
\begin{proof} We have $f(z,\nu,\varepsilon)= f_1 + f_2$, where $$f_1^2 -
f_2^2 = (1+\tau(\nu,\varepsilon)z)^2 - (1-z)(1-\tau^2 z)=
(2\tau(\nu,\varepsilon) + \tau^2 + 1)z + (\tau^2(\nu,\varepsilon)
- \tau^2)z^2,$$ and so
$$f(z,\nu,\varepsilon)=\frac{2\tau(\nu,\varepsilon) + \tau^2 +1
+(\tau^2(\nu,\varepsilon) - \tau^2)z}{(f_1(z) -
f_2(z))\,z^{-1}}\,.$$ The function $h=f_1 - f_2$ is analytic in
$\mathbb{D}$, $h(0) = 0$ and $$|h(z)|\leq 1 +
|\tau(\nu,\varepsilon)| + \left(2(1+\tau^2)\right)^{1/2}\,.$$ An
elementary analysis of the function $\tau(\nu,\varepsilon)$ shows
that
\begin{equation}\label{5.4} |\tau(\nu,\varepsilon)|\leq
|\tau|<1,\end{equation} i.e., $|h|\leq 4$. By Schwarz's Lemma $|h|\leq
4|z|$ for $|z|\leq 1$ and hence the second inequality in (\ref{5.2})
follows from (\ref{5.4}).

The below bound for the function $g$ comes out in a similar
manner: $$g=g_1 + g_2,\quad g_1^2 - g_2^2 = (1-\xi\,z)^2 -
(1-z)(1-\tau^2\,z)=(\xi^2 - \tau^2)\,z^2,$$ by Schwarz's Lemma
$|g_1 - g_2|\leq 4|z|^2$, which leads to (\ref{5.3}).\end{proof}

Let us next turn to the kernel $\tilde A$ (\ref{2.51}), (\ref{2.7}).
\begin{lemma}\label{lem5.2} The following inequalities take place
\begin{equation}\label{5.5}
|\tilde A(\lambda;n,m)|\leq
\frac{K_1(s,d)|w(\lambda)|}{|1-w^4(\lambda)|}\,d_m, \qquad
\lambda\in\mathbb{C},\quad n,m\in\mathbb{Z}_+,
\end{equation}
where
\begin{equation}\label{5.6}
K_1(s,d):=\frac{256\,d^3}{(d-s)^4}\,;\quad d_m=|b_m - b_m^0| +
\left|a_{m-1}^0 - \frac{a_{m-1}\,c_{m-1}}{a_{m-1}^0}\right|.
\end{equation}
\end{lemma}
\begin{proof} Consider the expression
\begin{equation}\label{5.7}
G(\lambda;n,m)\,w^{m-n} = \frac{\chi_n(\lambda)\hat\chi_m(\lambda)
-
\hat\chi_n(\lambda)\chi_m(\lambda)\,w^{2(m-n)}}{\langle\psi,\,\hat\psi\rangle}\,.
\end{equation}
For the denominator in (\ref{5.7}) we have
\begin{eqnarray*}\langle\psi,\,\hat\psi\rangle & =
& \frac{2a_2^0}{s_2}\left(u^2(\lambda) -
1\right)^{1/2}=\frac{2a_1^0\,a_2^0}{\lambda
-\nu}\left(u^2(\lambda) - 1\right)^{1/2} \\ \, & = & \frac{d^2 -
s^2}{4(\lambda-\nu)}\left(w^{-2}(\lambda) -
w^2(\lambda)\right).\end{eqnarray*}
 The functions $\chi_n,\
\hat\chi_n$ in the numerator are 2-periodic, and also
$$\chi_{2n}(\lambda)=\hat\chi_{2n}=1,\quad \chi_{2n+1}(\lambda)=
m(\lambda)\,w^{-1}(\lambda),\quad \hat\chi_{2n+1}(\lambda) =\hat
m(\lambda)\,w(\lambda),\quad n\in\mathbb{Z}_+.$$ Thereby we have 4
different expressions for the right hand side (\ref{5.7}):

\begin{equation}
  \label{5.8} G(\lambda;2p,2q)\,w^{2(q-p)} = D
 \,\frac{(\lambda - \nu)\,w^2\,\left(1 -
 w^{4(q-p)}\right)}{1 - w^4}\,, \end{equation}
  \begin{equation}\label{5.9} G(\lambda;2p,2q+1)\,w^{2(q-p)+1} =
D\,\frac{\hat m(\lambda)(\lambda - \nu)\,w^2 -
m(\lambda)(\lambda -
 \nu)\, w^{4(q-p)+3}}{1 - w^4}\,, \end{equation}
 \begin{equation} \label{5.10} G(\lambda;2p+1,2q)\,w^{2(q-p)-1} =
D\,\frac{ m(\lambda)(\lambda - \nu)\,w - \hat
m(\lambda)(\lambda -
 \nu)\, w^{4(q-p)+1}}{1 - w^4}\,, \end{equation}
 \begin{equation}\label{5.11} G(\lambda;2p+1,2q+1)\,w^{2(q-p)} =
D\,\frac{m(\lambda)\hat m(\lambda)(\lambda -
\nu)\,w^2\,\left(1 -
 w^{4(q-p)}\right)}{1 - w^4}
\end{equation}
with $$D=\frac{4}{d^2 -s^2}\,.$$

 We can estimate right hand sides in (\ref{5.8}) -(\ref{5.11}) by using
 the explicit formulae for $w$, $m$, $\hat m$. Assume first that
 $|\lambda|\geq d$. Then
\begin{eqnarray*}
 w^2(\lambda) &=& u(\lambda) - \left(u^2(\lambda) - 1\right)^{1/2} =
\frac{1}{u(\lambda) + \left(u^2(\lambda) - 1\right)^{1/2}} \\
   &=&\frac{d^2 -s^2}{2}\,\frac{1}{\lambda^2 - (s^2 + d^2 )/2 + \left((\lambda^2 -
s^2)(\lambda^2 - d^2)\right)^{1/2}}\,,
\end{eqnarray*}
so that
$$(\lambda - \nu)^2w^2(\lambda) = \frac{d^2 -
s^2}{2\,g(z)}\,\left(1 - \frac{\nu}{\lambda}\right)^2\,,\quad
z=\frac{d^2}{\lambda^2}\,.$$ By lemma \ref{lem5.1}
\begin{equation}\label{5.12}|(\lambda -
\nu)w(\lambda)|\leq\frac{\sqrt{2}\,d^2}{ \left(d^2 -
s^2\right)^{1/2}}\left(1 + \frac{|\nu|}{d}\right)\,,\quad
|\lambda|\geq d.\end{equation} The Weyl function is easily seen to
satisfy as $  z=\frac{d^2}{\lambda^2}$
\begin{eqnarray*}
  (\lambda - \nu)m(\lambda) &=& \frac{2a_1^0}{d^2 - s^2}\,\left(\lambda^2 +
d^2\,\tau(\nu,\varepsilon) - \left((\lambda^2 - s^2)(\lambda^2
-d^2)\right)^{1/2}\right)  \\
   &=&\frac{2a_1^0}{d^2 -s^2}\,\frac{2d^2\,\tau(\nu,\varepsilon) + s^2 + d^2 +
(d^4\,\tau^2(\nu,\varepsilon)-s^2\,d^2)\lambda^{-2}}{f(z,\nu,\varepsilon)},
.
\end{eqnarray*}
By the definition of $\tau(\nu,\varepsilon)$ (\ref{5.1}) and the property
(\ref{5.4}) it is not hard to derive
$$|2d^2\,\tau(\nu,\varepsilon) + s^2 +d^2|\leq (s+d)^2,\quad
|d^4\,\tau^2(\nu,\varepsilon) - s^2\,d^2|\leq 3\nu^2\,(s+d)^2,$$ and so
the application of lemma \ref{lem5.1} gives
\begin{equation}|m(\lambda)\,(\lambda - \nu)|\leq\frac{8a_1^0\,(s+d)^2}{d^2 -
s^2}\,\frac{4d^2}{(d-s)^2}\leq\frac{16(s+d)^2\,d^2}{(d-s)^3},\ \quad
|\lambda|\geq d \label{5.13}\end{equation} (in the later inequality we use
$2a_1^0\leq s+d$ (see (\ref{4.25})). In the same way
\begin{eqnarray*}
  \hat m(\lambda)(\lambda - \nu)w^2(\lambda) &=& \frac{2a_1^0}{d^2 -
s^2}\,\frac{\lambda^2 + d^2\,\tau(\nu,\varepsilon) +
\left(\lambda^2 - s^2)(\lambda^2 - d^2)\right)^{1/2}}{\lambda^2 -
d^2\,\xi + \left((\lambda^2 - s^2)(\lambda^2 -
d^2)\right)^{1/2}}\,\frac{d^2 - s^2}{2}
\\
   &=&a_1^0\,\frac{f(z,\nu,\varepsilon)}{g(z)}\,,\quad
   z=\frac{d^2}{\lambda^2}\,,
\end{eqnarray*}
and by lemma \ref{lem5.1}
\begin{equation}\label{5.14}
|\hat m (\lambda)(\lambda - \nu)w^2(\lambda)|\leq\frac{16a_1^0\,d^4}{(d^2
- s^2)^2}\leq\frac{8d^4}{(d-s)(d^2 - s^2)}\,,\quad |\lambda|\geq d.
\end{equation}
Finally
$$m(\lambda)\hat m(\lambda) =
-\frac{c_3(\lambda)}{s_2(\lambda)}=\frac{\lambda + \nu}{\lambda - \nu
},\quad m(\lambda)\hat m(\lambda)(\lambda - \nu)w(\lambda) = (\lambda +
\nu)w(\lambda)$$ and as above in (\ref{5.12})
\begin{equation}\label{5.15}
|m(\lambda)\hat m(\lambda)(\lambda -
\nu)w(\lambda)|\leq\frac{\sqrt{2}\,d^2}{ \left(d^2 -
s^2\right)^{1/2}}\left( 1 +\frac{|\nu|}{d}\right)\,,\quad
|\lambda|\geq d.
\end{equation}
For $|\lambda|\leq d$ the argument is much more elementary and is based on
$|w|\leq 1$:
\begin{equation}\label{5.16}
|(\lambda\pm \nu)w|\leq d+|\nu|,\end{equation}
\begin{equation}\label{5.17}
|m(\lambda)(\lambda - \nu)|\leq\frac{2a_1^0}{d^2 - s^2}\,\left(d^2
+d^2\,|\tau| + 2d^2\right)\leq\frac{4d^2}{d-s}\,,
\end{equation}
and the same estimate holds for $\hat m(\lambda)(\lambda - \nu)$.

It is clear that the worst bound comes in (\ref{5.13}) and thus we arrive
at the following inequalities
\begin{equation}\label{5.18}
|G(\lambda;n,m)w^{m-n}|\leq K_1(s,d)\,\frac{|w(\lambda)|}{|1 -
w^4(\lambda)|}\,,\quad \lambda\in\mathbb{C},\quad n,m\in\mathbb{Z}_+.
\end{equation}
The desired result now stems from (\ref{5.18}) and the definition
of the kernel $\tilde A$.\end{proof}

An application of (\ref{5.12}) and (\ref{4.25}) to (\ref{4.19}) gives
\begin{eqnarray}
  \nonumber |\chi_1(\lambda)| &\leq & \frac{2(d+s)}{d^2 -
 s^2}\,\frac{2\sqrt{2}\,d^2}{ \left(d^2 -
 s^2\right)^{1/2}}\,\frac{1}{|w^2(\lambda) - w^{2\varepsilon}(\nu)|} \\
  \label{5.19} &=& \frac{4\sqrt{2}\,d^2}
  {(d-s)\left(d^2 - s^2\right)^{1/2}}\,\frac{1}{|w^2(\lambda) -
 w^{2\varepsilon}(\nu)|}\,,
\end{eqnarray}
 and in fact this bound is true for all $\chi_n$.

 We are now in a position to estimate the solution of equation (\ref{2.6})
 \begin{equation}\label{5.20}
 V_n(\lambda)=\sum_{m=n+1}^\infty \tilde
 A(\lambda;n,m)\chi_m(\lambda) + \sum_{m=n+1}^\infty \tilde
 A(\lambda;n,m) V_m(\lambda)
 \end{equation}
 with constants expressed explicitly in terms of the spectral data for the
 2-periodic background. Assume first that
 \begin{equation}\label{5.21} \sum_{m=1}^\infty\, d_m < \infty.
 \end{equation}
 By (\ref{5.19}) and lemma \ref{lem5.2} the first series in (\ref{5.20})
converges uniformly on compact subsets of the domain
\begin{equation}\label{5.30}
 \Omega(\varepsilon)=\left\{\begin{array}{ll}
 \mathbb{C}\setminus\{\pm s,\pm d\},& \varepsilon = -1,\\
 \mathbb{C}\setminus\{\pm s,\pm d,\,\nu\},& \varepsilon =
 1,\end{array}\right.\end{equation}
 and its sum is subject to
 $$|V_{n,1}(\lambda)|\leq\frac{K_2(s,d)\kappa_0(n)|w(\lambda)|}{|w^2(\lambda)
 - w^{2\varepsilon}(\nu)||1-w^4(\lambda)|}\, ,$$
 \begin{equation}\label{5.211}
 K_2(s,d) = K_1(s,d)\,\frac{4\sqrt{2}\,d^2}{(d-s)\left(d^2 - s^2\right)^{1/2}},
 \end{equation}
 $\kappa_0(n)$ is defined in (\ref{2.8}), $\kappa(n) = o(1)$ as $n\to \infty$.
Hence for $V_n$ we have (cf.(\ref{2.10}))
  \begin{equation}\label{5.22}
|V_n(\lambda)|\leq \frac{K(\lambda)\kappa_0(n)}{|1-w^{4}(\lambda)|}
\exp\left(\frac{K(\lambda)\kappa_0(n)}{|1-w^{4}(\lambda)|}\right),\quad
K(\lambda) = \frac{K_2(s,d)|w(\lambda)|}{|w^2(\lambda) -
w^{2\varepsilon}(\nu)|},
\end{equation}
which leads to the following result.
\begin{theorem}\label{theor5.3}
In assumption $(\ref{5.21})$ there exists a solution $v_n$ of equation
$(\ref{2.4})$ which satisfies
\begin{equation}\label{5.23}
|v_n(\lambda) - \psi_n(\lambda)|\leq
\frac{K(\lambda)\kappa_0(n)|w^n(\lambda)|}{|1-w^{4}(\lambda)|}
\exp\left(\frac{K(\lambda)\kappa_0(n)}{|1-w^{4}(\lambda)|}\right),\quad
\lambda\in\Omega(\varepsilon),\end{equation} where $K(\lambda)$ is given
in $(\ref{5.22})$, $(\ref{5.211})$ and $(\ref{5.6})$, $\kappa_0(n)$ is
defined in $(\ref{2.8})$. \hfill $\square$\end{theorem}

Inequality (\ref{5.23}) with $n=0$ implies
$$|v_0(\lambda) - 1|\leq \frac{K(\lambda)}{|1 - w^4(\lambda)|}
\sum_{m=1}^\infty d_m\,\exp\left(\frac{K(\lambda)}{|1 - w^4(\lambda)|}
\sum_{m=1}^\infty d_m\right).$$

Taking into account the spectral interpretation of zeros of the Jost
function $v_0$ we can draw the following conclusion.
\begin{remark}\label{cor2}
Define the domain $G$ in $\lambda$-plain by inequality
\begin{equation}
G:=\left\{\lambda\in\mathbb{C}: \frac{K_2(s,d)|w(\lambda)|}{|w^2(\lambda)
-w^{2\varepsilon}(\nu)||1-w^4(\lambda)|}\,\sum_{m=1}^\infty\,d_m<t\right\},\label{5.31}
\end{equation}
where $t\simeq 0.567..$ is a unique positive root of the equation
$x\exp(x)=1$. Then $G$ is free from the discrete spectrum of $J$.
\hfill $\square$
\end{remark}

In the case when the first moment of the perturbation is finite
\begin{equation}
\sum_{m=1}^\infty \, m d_m<\infty, \label{5.24}\end{equation} the bounds
for the Green kernel are performed in a slightly different way by using a
simple inequality
$$ \left|\frac{1-w^n}{1-w}\right|\leq n,\quad |w|\leq 1.$$
For the right hand sides in (\ref{5.8}) and (\ref{5.11}) we see that
\begin{eqnarray*}
  |G(\lambda; 2p,2q)w^{2(q-p)}| &\leq& \frac{4(q-p)}{d^2 -s^2}\,|(\lambda -
\nu)w(\lambda)||w|, \\
  |G(\lambda; 2p +1,2q+1)w^{2(q-p)}| &\leq& \frac{4(q-p)}{d^2
-s^2}\,|(\lambda + \nu)w(\lambda)||w|.
\end{eqnarray*}
Next, for (\ref{5.9}) $$ G(\lambda;2p,2q+1)w^{2(q-p)+1}= G' +
G'',$$ where \begin{eqnarray*}G' & := & \frac{4(\lambda -
\nu)}{(d^2 - s^2)(w^{-2}(\lambda) - w^2(\lambda))}\,\left(\hat
m(\lambda) - m(\lambda)\right)w(\lambda)\\ \, & = &
\frac{s_2}{2a_2^0\left(u^2
-1\right)^{1/2}}\,\frac{\pm2\left(u^2-1\right)^{1/2}}{s_2}\,w(\lambda),
\end{eqnarray*}
and so
$$|G'|\leq\frac{|w(\lambda)|}{a_2^0}<\frac{2|w(\lambda)|}{d-s};$$
$$G'':=\frac{4(\lambda - \nu)}{d^2 - s^2}\,
\frac{m(\lambda)w(\lambda)\left(1 -
w^{4(p-q)}\right)}{w^{-2}(\lambda) - w^2(\lambda)},$$ and
$$|G''|\leq\frac{4(q-p)}{d^2 - s^2}\,\left|(\lambda -
\nu)mw^3\right|.$$ Hence $$|G(\lambda;2p,2q+1)\,w^{2(q-p) +1}|\leq
\frac{2|w|}{d-s} + \frac{4|(\lambda - \nu) m w^3|}{d^2 -
s^2}\,(q-p).$$ The similar inequality holds for
$G(\lambda;2p+1,2q)w^{2q - 2p - 1}$. Finally, $$|\tilde
A(\lambda;n,m)|\leq K_1(s,d)(m-n)d_m|w(\lambda)|.$$

The method of successive approximations of Section 3 is applicable and
provides the bound for the solution of the basic integral equation in
assumption (\ref{5.24}):
$$|v_n(\lambda) - \psi_n(\lambda)|\leq
K(\lambda)\kappa_1(n)\exp(K(\lambda)\kappa_1(n)),$$ where $K(\lambda)$ is
defined in (\ref{5.22}),
 $$\kappa_1(n):= \sum_{m=n+1}^\infty(m-n)d_m.$$

We are aimed here at the condition which guarantees the lack of discrete
spectrum for the operator $J$. That is why it seems reasonable to suppose
that the background operator $J_0$ has pure continuous spectrum, i.e., in
our notation  $\varepsilon = -1$. Then
 $$\frac{1}{|w^2(\lambda) -
 w^{-2}(\nu)|}=\frac{|w^2(\nu)|}{|1-w^2(\nu)w^2(\lambda)|}\leq\frac{1}{1-|w^2(\nu)|},$$
 and hence
 $$|v_0(\lambda) -
 1|\leq\frac { K_2(s,d)}{1-|w(\nu)|^2}\sum_{m=1}^\infty md_m\exp\left(\frac {K_2(s,d)
 }{1-|w(\nu)|^2}\sum_{m=1}^\infty md_m\right),$$
that gives the following result.
\begin{theorem}\label{theor5.4} The operator $J$ has no discrete spectrum as long as
$$\frac{K_2(s,d)}{1-|w(\nu)|^2}\sum_{m=1}^\infty md_m <t.$$
\end{theorem}

\end{document}